\documentclass[11pt]{amsart}
\makeindex
\usepackage{graphicx}
\usepackage[latin1]{inputenc} \usepackage[T1]{fontenc} \usepackage{lmodern}
\usepackage{pgf}
\usepackage{epsfig,graphicx,color,amsmath,a4wide,amssymb,amsfonts,amsbsy,xy,latexsym,epsf,pstricks,verbatim,times}

\usepackage{pgf,color}

\usepackage{changebar}

\definecolor{LightGrey}{rgb}{.85,.85,.85}
\definecolor{DarkGrey}{rgb}{.5,.5,.5}

\definecolor{Blue}{rgb}{.0,.0,0.9}
\definecolor{LightBlue1}{rgb}{.2,.4,0.9}
\definecolor{LightBlue2}{rgb}{.3,.5,0.9}
\definecolor{LightBlue3}{rgb}{.4,.6,0.9}
\definecolor{LightBlue4}{rgb}{.5,.7,.9}
\definecolor{LightBlue5}{rgb}{.6,.8,.9}
\definecolor{LightBlue6}{rgb}{.7,.9,.9}

\definecolor{Red}{rgb}{.9,.0,.0}
\definecolor{LightRed1}{rgb}{0.9,.2,.4}
\definecolor{LightRed2}{rgb}{0.9,.3,.5}
\definecolor{LightRed3}{rgb}{0.9,.4,.6}
\definecolor{LightRed4}{rgb}{.9,.5,.7}
\definecolor{LightRed5}{rgb}{.9,.6,.8}
\definecolor{LightRed6}{rgb}{.9,.7,.9}

\def\Jac{\mathop{\rm Jac }}

\definecolor{Grey}{rgb}{.5,.5,.5}

\definecolor{Blue}{rgb}{.0,.0,0.9}
\definecolor{LightBlue1}{rgb}{.2,.4,0.9}
\definecolor{LightBlue2}{rgb}{.3,.5,0.9}
\definecolor{LightBlue3}{rgb}{.4,.6,0.9}
\definecolor{LightBlue4}{rgb}{.5,.7,.9}
\definecolor{LightBlue5}{rgb}{.6,.8,.9}
\definecolor{LightBlue6}{rgb}{.7,.9,.9}

\definecolor{Red}{rgb}{.9,.0,.0}
\definecolor{LightRed1}{rgb}{0.9,.2,.4}
\definecolor{LightRed2}{rgb}{0.9,.3,.5}
\definecolor{LightRed3}{rgb}{0.9,.4,.6}
\definecolor{LightRed4}{rgb}{.9,.5,.7}
\definecolor{LightRed5}{rgb}{.9,.6,.8}
\definecolor{LightRed6}{rgb}{.9,.7,.9}

\xyoption{all}
\usepackage[english]{babel}

\newcounter{noalgo}[section]
\newdimen\indentalgo
\newdimen\indentalgodec\indentalgo=0.0mm\indentalgodec=10mm

\newcommand{\If}{\advance\indentalgo by \indentalgodec {\bf if }}

\newcommand{\For}{\global\advance\indentalgo by \indentalgodec {\bf for }}

\newcommand{\Endindent}{\global\advance\indentalgo by -\indentalgodec}

\newdimen\decalage \decalage=0.5cm
\newcounter{algo} 

\newcommand{\PP}{\mathbb P}
\newcommand{\bb}{\mathbf b}

\def\<<{\leavevmode
  \raise0.28ex\hbox{$\scriptscriptstyle\langle\!\langle$}\nobreak
  \hskip -.6pt plus.3pt minus.2pt\,}
\def\>>{\,\nobreak\hskip -.6pt plus.3pt minus.2pt
  \raise0.28ex\hbox{$\scriptscriptstyle\rangle\!\rangle$}}

\def\Aut{\mathop{\rm{Aut}}\nolimits }

\def\Hom{\mathop{\rm{Hom}}\nolimits }
\def\Roots{\mathop{\rm{Roots}}\nolimits }
\def\Bij{\mathop{\rm{Bij}}\nolimits }
\def\Alg{\mathop{\rm{Alg}}\nolimits }
\def\Sym{\mathop{\rm{Sym}}\nolimits }

\def\PGL{\mathop{\rm{PGL}}\nolimits }
\def\GL{\mathop{\rm{GL}}\nolimits }

\def\GL{\mathop{\rm{GL}}\nolimits }

\def\FF{{\mathbb F}}

\def\Fq{{\FF _q}}

\def\Fqs{{{\FF ^*_q}}}

\def\bm{{\bf m  }}
\def\bx{{\bf x  }}
\def\bX{{\bf  X  }}
\def\by{{\bf y  }}

\def\QQ{{\mathbb Q}}

\def\RR{{\mathbb R}}

\def\cI{{\mathcal I}}

\def\cO{{\mathcal O}}
\def\cP{{\mathcal P}}

\def\cS{{\mathcal S}}

\newtheorem{lemma}{Lemma}

\newtheorem{theorem}{Theorem}

\providecommand{\myproofname}{Proof}

\usepackage{enumitem}

\begin{document}
\begin{abstract}
  We explore  parameterizations by radicals  of low genera algebraic curves.
  We prove that for $q$ a prime power that is large enough and prime to $6$, a
  fixed positive proportion of all genus 2 curves over the field with $q$
  elements can be parameterized by $3$-radicals.
  This results in the existence of a deterministic encoding into these curves
  when $q$ is congruent to $2$ modulo $3$.
  We extend this construction to parameterizations by $\ell$-radicals for small odd
  integers $\ell$, and make it explicit for $\ell=5$.
\end{abstract}

\title{The geometry of some parameterizations and encodings}

\author{Jean-Marc Couveignes}
\address{Jean-Marc Couveignes, Univ. Bordeaux, IMB, UMR 5251, F-33400 Talence, France.}
\address{Jean-Marc Couveignes, CNRS, IMB, UMR 5251, F-33400 Talence, France.}
\address{Jean-Marc Couveignes, INRIA, F-33400 Talence, France.}
\address{Jean-Marc Couveignes, Laboratoire International de Recherche en Informatique et Math\'ematiques Appliqu\'ees.}
\email{Jean-Marc.Couveignes@math.u-bordeaux1.fr}

\author{Reynald Lercier}

\address{%
  \textsc{Reynald Lercier, DGA MI}, %
  La Roche Marguerite, %
  35174 Bruz, %
  France. %
}
\address{Reynald Lercier, 
  Institut de recherche math\'ematique de Rennes, %
  Universit\'e de Rennes 1, %
  Campus de Beaulieu, %
  35042 Rennes, %
  France. %
} 
\email{reynald.lercier@m4x.org}

\thanks{Research supported by the ``Direction  G{\'e}n{\'e}rale de 
    l'Armement'', by the ``Agence Nationale de la Recherche'' (project PEACE),
and  the Investments for the future Programme IdEx Bordeaux (ANR-10-IDEX-03-02)
through the CPU cluster (Numerical certification and reliability).} 

\date{\today}

\maketitle

\section{Introduction}

Let $\FF_{q}$ be a finite field, let $C/{\FF_{q}}$ be an algebraic curve, we
propose in this paper new algorithms for computing in deterministic polynomial
time a point in $C(\FF_{q})$. This is useful in numerous situations, for
instance in discrete logarithm cryptography~\cite{BonFra2001}.
To be more precise, we consider this question for low genus curves with an
emphasis on the genus 2 case.
The mathematical underlying problem is to compute radical expressions for
solutions of a system of algebraic equations.
Galois theory  provides nice answers, both in theory and practice,
for sets of dimension 0 and   degree less than 5.
Explicit results are known in dimension 1 too.
A famous theorem of Zariski states that a generic curve of genus at least 7
cannot be parameterized by radicals.
Conversely, a complex curve of genus less than 7 can be parameterized by
radicals over the field of rational fractions \cite{Harrison,Fried}.

In this work we restrict the degree of radicals involved in the
parameterizations.
Typically, for $C$ a curve over the field with $q$ elements, we only allow
radicals of degrees $l$ prime to $q(q-1)$.
The reason is that for such $l$, we can compute $l$-th roots of elements in
$\FF_q$ in deterministic polynomial time in $\log(q)$.
Especially,  we do not allow square roots.  We will be mainly
concerned with genus 2 curves.

%\medskip

Following pioneering investigations \cite{SK,SKS} by Schinzel and  Ska{\l}ba,
Shallue
and Woestijne  came in 2006  to a first practical deterministic algorithm
for constructing points on genus 1 curves over any finite field~\cite{SW}.
In 2009, Icart proposed  a deterministic encoding with quasi-quadratic
complexity in $\log q$ for elliptic curves over a finite field when $q$ is
congruent to $2$ modulo $3$~\cite{icart}. To this end, he  constructed
a parameterization by $3$-radicals for
every elliptic curve over a field with characteristic prime  to $6$.
Couveignes and Kammerer recently proved that there exists an infinity of
such parameterizations~\cite{CK12}, corresponding to rational curves on a K3
surface associated with the elliptic curve.
Nevertheless, in genus 2, only partial results are known.
Ulas  attempted   to generalize Shallue and Woestijne
results~\cite{ulas07}.
Tibouchi and Fouque designed encodings for curves with automorphism group
containing the dihedral group with 8 elements~\cite{FM}.
Each of these two constructions reaches a family of dimension 1 inside the
dimension 3 moduli space of the genus 2 curves. So the proportion of target
curves for such parameterizations tends to zero when $q$ tends to infinity.
At the same time, Kammerer, Lercier and Renault~\cite{KLR} published encodings
for a dimension 2 family of genus 2 curves.
In particular their curves have no non-hyperelliptic involution.  However
these curves still represent a negligible proportion of all genus 2 curves
when $q$ tends to infinity.

In this paper we construct a parameterization by $3$-radicals for a genus 2
curve $C$ over a field $K$ with characteristic $p$ prime to $6$ under the
sole restriction that $C$ has two $K$-rational points whose difference has
order $3$ in the Jacobian variety.  This is a dimension 3 family.
In particular, we parameterize all genus 2 curves when $K$ is
algebraically closed. When $K$ is a finite field with characteristic
prime to $6$ we parameterize
a positive proportion of all curves in that way.
Our construction extends the ones by Farashahi \cite{fara} for genus
1 curves
and  Kammerer, Lercier and Renault~\cite{KLR} for genus 2 curves.
Our starting point is the observation that the role played by Tartaglia-Cardan
formulae in these  parameterizations can be formalized and generalized
using the theory of torsors under solvable finite group schemes.
This leads us to a systematic exploration and combination of possibilities
offered by the action of small solvable group schemes over curves of low
genus.
%
% For example, we prove that for $q$ a prime power that is large enough and
% prime to $6$, a fixed positive proportion of all genus 2 curves over the
% field with $q$ elements can be parameterized by $3$-radicals.
%
% This results in the existence of a deterministic encoding into these curves
% when $q$ is congruent to $2$ modulo $3$.
%

The principles of our method are presented in Section~\ref{sec:defgene}.
We first recall the basics of parameterizations of curves by radicals and
encodings, then we explain how to produce such parameterizations using the
action of solvable group schemes on algebraic curves.
Section~\ref{sec:g1} provides a first illustration of this general method in
the case of genus 1 curves. This offers a new insight on previous work by
Farashahi, and Kammerer, Lercier, Renault.
We present a parameterization of a 3~~dimensional family of genus
2 curves in Section~\ref{sec:g2}.  Variations on this theme are presented in
Section~\ref{sec:other}.
Section~\ref{sec:g25} presents detailed computations for one of these
families.  We parameterize Jacobians of dimension 2 with one 5-torsion
point.  We finish with a few questions and prospects.
%

%\medskip

We thank Jean Gillibert, Qinq Liu, and Jilong Tong for useful discussions.

\section{Definitions and generalities}\label{sec:defgene}

In this section we recall a few definitions and present the principles
of our method. Sections~\ref{sec:radic}  and~\ref{sec:radic2}
recall elementary results about radicals. Section~\ref{sec:paramet}
recalls the definition of a parameterization. Section~\ref{sec:torsors}
gives elementary definitions about torsors. Basic properties
of encodings are recalled in Section~\ref{sec:encodings}.
Section~\ref{sec:cardan} presents Tartaglia-Cardan formulae in the
natural language of torsors. Our strategy for finding new parameterizations
is presented in Sections~\ref{sec:mu3mu2} and~\ref{sec:selec}.

\subsection{Radical extensions}\label{sec:radic}

The following classical lemma  \cite[Chapter VI, Theorem 9.1]{lang}
gives necessary and sufficient conditions for a binomial to be irreducible.

\begin{lemma}\label{lem:cond}
Let $K$ be a field,  let $d\ge 1$ be a
positive integer,  and let $a\in K^*$. The  polynomial $x^d-a$
is irreducible in $K[x]$
if and only if the two following conditions hold true
\begin{itemize}
\item For every prime integer $l$ dividing $d$, the scalar $a$ is not the  $l$-th power of an element in $K^*$,
\item If $4$ divides $d$, then $-4a$ is not the  $4$-th power of an element in  $K^*$.
\end{itemize}
\end{lemma}

Let $K$ be a field with characteristic $p$. Let $S$ be a
set of rational primes such that $p\not\in S$.
Let $M\supset K$
be a  finite separable $K$-algebra,  and $L\subset M$ a $K$-subalgebra
of $M$.
The extension $L\subset M$ is said to be
$S$-{\it radical} if $M$ is isomorphic, as an $L$-algebra,
to $L[x]/(x^l-a)$ for some
$l \in S$
and some $a\in L^*$.
When  $S$ contains all  primes but $p$, we speak of {\it radical extensions}.

An  extension $M\supset L$ is said to be
$S$-{\it multiradical}
if  there exists a finite sequence of $K$-algebras
\[K\subset L=L_0\subset L_1  \subset \dots \subset L_n=M\]
 such that every intermediate extension $L_{i+1}/L_i$
for $0\le i\le n-1$ is $S$-radical.

\subsection{Radical morphisms}\label{sec:radic2}

Let $K$ be a field with characteristic $p$.
Let $\bar K\supset K$ be an algebraic closure.
Let $f : C \rightarrow D$ be an epimorphism of
(projective, smooth, absolutely integral) curves
over $K$.
We say that $f$ is a {\it radical morphism} if
the associated function field extension
$K(D)\subset K(C)$
is radical. We define similarly multiradical morphisms,
$S$-radical morphisms,
$S$-multiradical morphisms.
If $f$ is a radical  morphism then $K(C)=K(D,b)$
where $b^l=a$ and
$a$ is a non-constant function on $D$ and $l\not = p$ is a prime
integer.
%Let $\Delta : C\rightarrow C\times C$ the diagonal morphism.
%Call $\gamma_b$ the composition of $\Delta$ and $f\times b$
Call $\gamma_b$ the map
\begin{equation*}\label{eq:gamma}
    \xymatrixrowsep{0.1in}
    \xymatrixcolsep{0.3in}
\xymatrix{
\gamma_b : & C \ar@{->}[r]    & D\times \PP^1\\
&P\ar@{|->}[r]&(f(P),b(P)).
}
\end{equation*}
Let $X\subset C$ be the ramification  locus of $f$, and  let $Y=f(X)\subset D$  be the branch locus.
A geometric point $Q$ on $D$ is branched if and only if
$a$ has a zero or a pole at $Q$ with multiplicity  prime to $l$.
We ask if $\gamma_b$
induces an  injection on
$C(\bar K)$.
Equivalently we ask if $b$
separates points in every fiber
of $f$.
First, there is a unique ramification point
above each branched point. Then, if $a$ has neither a  zero nor a pole
at $Q$, then  $b$ separates the points in the fiber of $f$ above $Q$.
Finally, if $a$ has a zero or a pole at $Q$ with multiplicity
divisible by $l$, then $b$ (and $\gamma_b$) fail to separate
the points in the fiber of $f$ above $Q$.
However, there exists a finite covering $(U_i)_i$ of $C$
by affine open subsets,  and functions $b_i \in \cO(U_i  -  X)^*$
such that $b_i/b\in K(D)^*\subset K(C)^*$. We set
$\bb=(b_i)_{1\le i\le I}$ and define a map
\begin{equation*}\label{eq:gammab}
    \xymatrixrowsep{0.1in}
    \xymatrixcolsep{0.3in}
\xymatrix{
\gamma_\bb : & C \ar@{->}[r]    & D\times \left( \PP^1\right)^I\\
&P\ar@{|->}[r]&(f(P),b_1(P), \dots, b_I(P)).
}
\end{equation*}
This map induces  an injection on $C(\bar K)$. So every point $P\in C(\bar K)$ can
be characterized by its image $f(P)$ on $D$ and the value of the $b_i$ at $P$.
%Indeed, one such $b_i(P)$ suffices (take one $i$ such that $a_i(Q)$ is finite and non-zero).

\subsection{Parameterizations}\label{sec:paramet}

An $S$-{\it parameterization}   of a projective, absolutely integral,
smooth curve $C$  over $K$ is a triple  $(D,\rho,\pi)$
where $D$ is  another
projective, absolutely integral,
smooth curve  over $K$, and $\rho$ is  an $S$-multiradical map
from $D/K$ onto $\PP^1/K$,  and $\pi$ is an  epimorphism from
$D/K$  onto $C/K$.
In this situation one says that $C/K$ is {\it parameterizable} by $S$-radicals.

\begin{equation}\label{eq:para}
\xymatrix{
&D \ar@{->}[dl]_\pi\ar@{->}[d]^\rho    \\
C&\PP^1
}\end{equation}

\subsection{$\Gamma$-groups}\label{sec:torsors}

Let $K$ be a field with characteristic $p$.
Let  $K_s$ be  a separable  closure of $K$. Let $\Gamma$ be the Galois
group of $K_s/K$. Let $A$ be a finite set acted on
continuously
by $\Gamma$. We say that
 $A$ is  a finite  $\Gamma$-set.
We associate to it the separable $K$-algebra
\[\Alg(A)=\Hom_\Gamma(A,K_s)\]
of $\Gamma$-equivariant maps from $A$ to $K_s$.
If $G$ is a finite $\Gamma$-set and
has a group structure compatible with  the $\Gamma$-action
we say that $G$ is a finite $\Gamma$-group, or a finite {\'e}tale
group scheme over
$K$. Now let  $A$ be a finite  $\Gamma$-set acted on
 by a  finite $\Gamma$-group $G$. If the
action of $G$ on $A$ is compatible with the actions of $\Gamma$ on $G$ and $A$, then we say
that $A$ is a  finite $G$-set.
The quotient $A/G$ is then a finite $\Gamma$-set.
If further $G$ acts freely
on $A$ we say that  $A$ is a  free  finite $G$-set.
A simply transitive $G$-set is called a $G$-torsor.
The left action of $G$ on itself defines a $G$-torsor called
the trivial torsor. The set of isomorphism classes of $G$-torsors is isomorphic, as a pointed
set, to $H^1(\Gamma,G)$. See \cite[Chapter I \S 2]{NSW}.

Let  $l\not=p$ be a prime and let  $A$ be a free  finite $\mu_l$-set.
Let $B=A/\mu_l$. According to Kummer theory,
the inclusion  $\Alg(B)\subset\Alg(A)$ is a radical
extension of separable $K$-algebras. It has degree $l$.

Let $S$ be a finite set of primes.
Assume that the characteristic $p$ of $K$
does not belong to $S$.
A finite $\Gamma$-group $G$  is said to be $S$-{\it solvable} if there exists
a sequence of $\Gamma$-subgroups
$1=G_0\subset G_1\subset \dots  \subset G_I=G$ such that for every $i$ such that $0\le i\le I-1$, the
group $G_i$ is normal in $G_{i+1}$, and the
quotient $G_{i+1}/G_i$ is isomorphic, as a finite $\Gamma$-group,  to $\mu_{l_i}$
for some $l_i$ in $S$.

 Let
$G$ be a finite $\Gamma$-group. Assume that  $G$  is  $S$-{\it solvable}.
Let  $A$ be a free finite $G$-set.
Let $B=A/G$. The inclusion  $\Alg(B)\subset\Alg(A)$ is an $S$-multiradical
extension of separable $K$-algebras. It has degree $\# G$.

\subsection{Encodings}\label{sec:encodings}

We assume that  $K$ is  a finite field with characteristic $p$ and cardinality $q$.
Let $S$ be a  set of prime integers. We assume that $p\not \in S$
and $S$
is disjoint from the support of $q-1$.
Let $C$ and $D$ be two projective, smooth, absolutely integral curves over $K$. Let $f : C\rightarrow D$ be a radical morphism of degree $l\in S$.
Let $X\subset C$ be the ramification  locus of $f$, and  let $Y=f(X)\subset D$  be the branch locus.
Let $F : C(K) \rightarrow D(K)$
the induced map on $K$-rational points.
We prove that $F$ is a bijection.

A branched point $Q$ in $D(K)$ is totally ramified, so has
 a unique preimage $P$ in $C(K)$. Let $Q\in D(K)-Y(K)$ be a non-branched point. The fiber $f^{(-1)}(Q)$
is a $\mu_l$-torsor. Since $H^1(K,\mu_l)=K^*/(K^*)^l$ is trivial, this torsor is
isomorphic to $\mu_l$ with the left action. Since
$H^0(K,\mu_l)=\mu_l(K)$ is trivial also,   $f^{(-1)}(Q)$ contains a unique $K$-rational point.
Therefore $F$ is a bijection.

\begin{lemma}\label{lem:bij}
Let $K$ be a finite field with $q$ elements. Let
$S$ be a finite set of prime integers. We assume that $p\not \in S$
and $S$ is disjoint from the support of $q-1$.
Let $f : C\rightarrow D$ be an $S$-multiradical morphism between
two smooth, projective, absolutely irreducible curves over $K$.
The induced map $F : C(K)\rightarrow D(K)$
on $K$-rational points is a bijection.
\end{lemma}

The reciprocal map $F^{(-1)} : D(K)\rightarrow C(K)$
can be evaluated in deterministic polynomial time by
computing successive $l$-th roots for various $l\in S$.

We assume now that we  are in the situation of the diagram  (\ref{eq:para}). Let
$R : D(K)\rightarrow \PP^1(K)$ be the map induced by $\rho$ and let
$\Pi : D(K)\rightarrow C(K)$ be the map induced by $\pi$.
The composition $\Pi \circ R^{(-1)}$ is called an {\it encoding}.

\subsection{Tartaglia-Cardan formulae}\label{sec:cardan}

Let $K$ be a field with characteristic prime to $6$.
Let $K_s$ be an algebraic closure of $K$.
Let $\Gamma$ be the Galois group of $K_s/K$.
Let $\mu_3 \subset  K_s$ be the finite $\Gamma$-set consisting of the three roots
of unity. Let $\Sym(\mu_3)$ be the full permutation  group on  $\mu_3$.
The Galois group  $\Gamma$ acts on $\mu_3$. So we have a group homomorphism
$\Gamma \rightarrow \Sym(\mu_3)$ and
$\Gamma$ acts on $\Sym(\mu_3)$ by conjugation. This action turns
$\Sym(\mu_3)$ into a group scheme over $K$. Because $\mu_3$ acts on itself
by translation, we have an  inclusion of group schemes
$\mu_3\subset \Sym(\mu_3)$ and $\mu_3$ is a normal subgroup of $\Sym(\mu_3)$.
The stabilizer of $1\in \mu_3$ is a subgroup scheme of  $\Sym(\mu_3)$. It is not
normal in  $\Sym(\mu_3)$. It is isomorphic to $\mu_2$. So   $\Sym(\mu_3)$ is
the semidirect product $\mu_3\rtimes \mu_2$.
Let $\zeta_3\in  K_s$ be a primitive third root of unity.
We set $\sqrt{-3}=2\zeta_3+1$.
Let \[h(x)=x^3-s_1x^2+s_2x-s_3\] be a degree $3$ separable polynomial
in $K[x]$. Let \[R=\Roots(h)\]  be the set of roots
of  $h(x)$ in
$K_s$.  This is a finite $\Gamma$-set with cardinality $3$.
Let \[A=\Bij(\Roots(h), \mu_3)\] be the set of bijections
from $R$ to $\mu_3$.
For $\gamma\in \Gamma$ and $f\in A$ we set ${}^\gamma f=
\gamma  \circ f \circ \gamma^{-1}$.
This turns $A$ into a  finite $\Gamma$-set of cardinality $6$.
The action of $\Sym(\mu_3)$ on the left turns it into a
$\Sym(\mu_3)$-torsor.
Let \[C=A/\mu_3\] be the quotient of $A$ by the normal $\Gamma$-subgroup
$\mu_3\subset \Sym(\mu_3)$ of order $3$. This is a $\mu_2$-torsor.
Let \[B=A/\mu_2\] be the  quotient of $A$ by the stabilizer
of $1$ in $\Sym(\mu_3)$. This is a finite $\Gamma$-set of cardinality~$3$,
naturally isomorphic to $\Roots(h)$.
We define  a function $\xi$ in $\Alg(B)\subset \Alg(A)$   by
\begin{equation*}\label{eq:xi}
    \xymatrixrowsep{0.1in}
    \xymatrixcolsep{0.3in}
\xymatrix{
\xi : & A\ar@{->}[r] & K_s\\
&f\ar@{|->}[r]&f^{(-1)}(1).
}\end{equation*}
The algebra $\Alg(B)$ is generated by $\xi$,  and
the characteristic polynomial of $\xi$ is $h(x)$.
So \[\Alg(B)\simeq K[x]/h(x).\]
Tartaglia-Cardan formulae construct functions in the algebra
$\Alg(A)$ of the $\Sym(\mu_3)$~-~torsor $A$.
These functions can be constructed with radicals
because \[\Sym(\mu_3)=\mu_3\rtimes \mu_2\] is $\{2,3\}$-solvable.
A  first function $\delta$ in $\Alg(C)\subset \Alg(A)$ is defined by
\begin{equation*}\label{eq:delta}
    \xymatrixrowsep{0.1in}
    \xymatrixcolsep{0.3in}
\xymatrix{
\delta : & A\ar@{->}[r] & K_s\\
&\scriptstyle f\ar@{|->}[r]&\scriptstyle \sqrt{-3}\left(f^{(-1)}(\zeta)-f^{(-1)}(1)\right)
\left(f^{(-1)}(\zeta^2)-f^{(-1)}(\zeta)\right)\left(f^{(-1)}(1)-f^{(-1)}(\zeta^2)\right).
}\end{equation*}
Note that the $\sqrt{-3}$ is necessary to balance the Galois action on $\mu_3$.
The algebra $\Alg(C)$ is generated by $\delta$. And
\[\delta^2 = 81s_3^2-54s_3s_1s_2-3s_1^2s_2^2+12s_1^3s_3+12s_2^3=-3\Delta\]
is the discriminant $\Delta$ of $h(x)$ multiplied by $-3$.
We say that $-3\Delta$ is the {\it twisted discriminant}.
A  natural function $\rho$ in $\Alg(A)$ is defined as
\begin{equation*}\label{eq:rho}
    \xymatrixrowsep{0.1in}
    \xymatrixcolsep{0.3in}
\xymatrix{
\rho : & A\ar@{->}[r] & K_s\\
&f\ar@{->}[r]&\sum_{r\in R}r\times f(r)=\sum_{\zeta \in \mu_3}\zeta
\times f^{(-1)}(\zeta).
}\end{equation*}
It is clear that $\rho^3$ is invariant
by $\mu_3\subset \Sym(\mu_3)$ or equivalently
belongs to $\Alg(C)$. So it can be expressed
as a combination of $1$ and $\delta$. Indeed a simple calculation shows that
\[\rho^3= s_1^3+\frac{27}{2}s_3-\frac{9}{2}s_1s_2-\frac{3}{2}\delta.\]
A variant of $\rho$ is
\begin{equation*}\label{eq:rho'}
    \xymatrixrowsep{0.1in}
    \xymatrixcolsep{0.3in}
\xymatrix{
\rho' : & A\ar@{->}[r] & K_s\\
&f\ar@{->}[r]&\sum_{r\in R}r\times f(r)^{-1}.
}\end{equation*}
One has
\[\rho'^3= s_1^3+\frac{27}{2}s_3-\frac{9}{2}s_1s_2+\frac{3}{2}\delta\]
and \[\rho\rho'=s_1^2-3s_2.\]
Finally, the root $\xi$ of $h(x)$  can be expressed
in terms of $\rho$ and $\rho'$ as
\begin{equation*}
\xi=\frac{s_1+\rho+\rho'}{3}.
\end{equation*}
Note that the algebra $\Alg(A)$ is not the Galois closure of $K[x]/h(x)$.
If we wanted to construct a Galois closure we would rather consider
the $\Sym(\{1,2,3\})$-torsor $\Bij(R,\{1,2,3\})$ of indexations of the roots. We are not
interested in this torsor however. This
is  because   $\mu_3\rtimes \mu_2$  is solvable while
$C_3\rtimes C_2$ is not, in general.
The algebra constructed by Tartaglia and Cardan contains the initial cubic
extension, because the quotient of $\Bij(\Roots(h), \mu_3)$ by the
stabilizer of $1$ in $\Sym(\mu_3)$
is isomorphic to  the quotient of $\Bij(R,\{1,2,3\})$ by the
stabilizer of $1$ in $\Sym(\{1,2,3\})$, that is
$\Roots(h)$. On the other hand, the quotient
 of $\Bij(R,\{1,2,3\})$ by the $3$-cycle $(123)\in \Sym(\{1,2,3\})$ is associated with the algebra
$K[x]/(x^2-\Delta )$ while the quotient
 of $\Bij(R,\mu_3)$ by the $3$-cycle $(1\zeta\zeta^2)\in \Sym(\mu_3)$ is associated with the algebra
$K[x]/(x^2+3\Delta )$.

\subsection{Curves with a $\mu_3\rtimes \mu_2$ action}\label{sec:mu3mu2}

We still assume that the characteristic of $K$ is prime to $6$.
Let $A$ be a projective, absolutely integral,
smooth curve  over $K$. We assume that
the automorphism group
 $\Aut(A\otimes_K K_s)$  contains
a finite {\'e}tale $K$-group-scheme isomorphic to
 $\mu_3\rtimes \mu_2$.
The quotients $B=A/\mu_2$, and $C=A/\mu_3$ are
 projective, absolutely integral,
smooth curves  over $K$.
In this situation, we say that $C$ is the {\it resolvent} of $B$.
By abuse of language we may say also that we have constructed
a parameterization of $B$ by $C$.

Assume now that  $C$ admits a  parameterization  by $S$-radicals as in diagram~(\ref{eq:para}).
We call $D'$ the normalization of the fiber product of $A$ and $D$ above $C$. We assume that
$D'$ is absolutely integral.
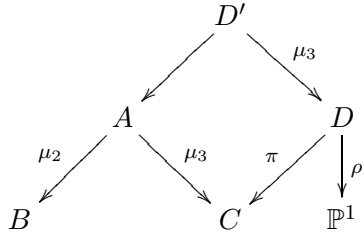
\begin{equation*}\label{eq:passage}
\xymatrix{
&&D' \ar@{->}[dl]\ar@{->}[dr]^{\mu_3} &   \\
&A\ar@{->}[dl]_{\mu_2}\ar@{->}[dr]^{\mu_3}&&D\ar@{->}[dl]_{\pi}\ar@{->}[d]^{\rho}\\
B&&C&\PP^1
}\end{equation*}
We set $S'=S\cup \{3\}$. We let $\rho'$ be the composite map
\[\rho' : D'\stackrel{\mu_3}{\longrightarrow} D\stackrel{\rho}{\longrightarrow} \PP^1,\]
and $\pi'$ the
 composite map
\[\pi' : D'\longrightarrow A\stackrel{\mu_2}{\longrightarrow} B.\]
Then $(D',\rho',\pi')$ is an $S'$-parameterization of $B$.
The mild condition that $D'$ be absolutely integral is granted in the following cases:
\begin{enumerate}
\item When $C=\PP^1$ and $\pi$ and $\rho$ are trivial.
\item When the $\mu_3$-quotient  $A\rightarrow C$ is branched at some point $P$ of $C$,  and $\pi$  is not branched at $P$. Indeed the two coverings  are linearly disjoint in that case. We
note that when
$C$ has genus 1 we may   compose $\pi$ with a translation
to ensure that it is not branched at $P$.
\item When the degree of $\pi$ is prime to $3$, because  $A\rightarrow C$ and
$\pi$ are linearly disjoint then. Note that the resulting parameterization
$\pi'$  has degree prime to $3$ also. We can iterate in that case.
\end{enumerate}

\subsection{Selecting curves}\label{sec:selec}

We still assume that the characteristic of $K$ is prime to~$6$.
We now look for interesting examples of curves with a $\mu_3\rtimes \mu_2$ action.
We keep the notation introduced in Section~\ref{sec:mu3mu2}. We set
$E=A/(\mu_3\rtimes \mu_2)$.
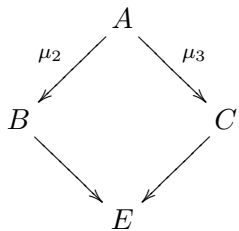
\begin{equation*}\label{eq:passage2}
\xymatrix{
&A\ar@{->}[dl]_{\mu_2}\ar@{->}[dr]^{\mu_3}&\\
B\ar@{->}[dr]&&C\ar@{->}[dl]\\
&E&
}\end{equation*}
The curve $C$ is the one we already know how to parameterize.
The curve $B$ is the one we want to parameterize. It should be as generic as possible.
In particular, we will assume that $E=\PP^1$. Otherwise, the Jacobian of $B$ would contain
a subvariety isogenous to the Jacobian of $E$. It would not be so generic then.

Assuming now that $E=\PP^1$ we denote by $r$ the number of branched points of the cover $B\rightarrow E$.
Let $r_s$ be the number of  branched points with ramification type $2,1$. These are called simple branched
points. Let $r_t$ the number  of  branched points with ramification type $3$. These are totally branched points.
We have $r=r_s+r_t$. According to the Hurwitz Genus Formula \cite[III.4.12, III.5.1]{Stich} the genus of $B$ is
\[g_B=\frac{r_s}{2}+r_t-2.\]
We note that every simple branched point of the cover
$B\rightarrow E$ gives rise to a branched point
of type $2,2,2$ of the cover $A\rightarrow E$ and to a (necessarily simple) branched point of $C\rightarrow E$.
And every totally  branched point of the cover $B\rightarrow E$ gives rise to a branched point
of type $3,3$ of the cover $A\rightarrow E$ and to a non-branched point of $C\rightarrow E$. So
\[g_A=\frac{3r_s}{2}+2r_t-5, \text{ \, and \,\,\,}g_C=\frac{r_s}{2}-1.\]
We set \[m=r-3=r_s+r_t-3\] and call it the {\it modular dimension}. It is the dimension of the family of covers  obtained
by  letting the $r$ branched points move along $E=\PP^1$. The $-3$ stands for the action of $\Aut(\PP^1)=\PGL_2$.
If we aim at all curves of genus $g_B$ we should have $m$ greater than or equal to the
dimension  of the moduli space of curves of genus $g_B$.
We deduce the {\it genericity condition}
\[r_s+4r_t\le 12-2\epsilon(\frac{r_s}{2}+r_t-2),\]
where $\epsilon(0)=3$, $\epsilon(1)=1$, and $\epsilon(n)=0$ for $n\ge 2$.
This is a necessary condition.

The first case to consider is when $C$ has genus 0 (because we know how to parameterize genus 0 curves). So we first take $r_s=2$.
So $g_B=r_t-1$ and the genericity condition reads $r_t\le 2$. Only $r_t=2$ is of interest. We shall see in Section~\ref{sec:g1}
that we find a  parameterization similar to those by
Farashahi and Kammerer, Lercier, Renault
in this case.

Assuming that we know how to parameterize some genus 1 curves, we may consider the case when $C$ itself has genus 1.
We have $r_s=4$ in that case. And $g_B=r_t$. The genericity assumption reads
$r_t\le 2$. The case $r_t=2$
will be studied in detail in Section~\ref{sec:g2}.

%JMC : Mais au fait, on n'a pas pens{\'e} {\`a} regarder le cas $b=1$ qui peut {\^e}tre int{\'e}ressant comme alternative {\`a} Icart.

\section{Curves of genus 1}\label{sec:g1}

Let $K$ be a field of  characteristic prime to $6$.
Let $B/K$ be a projective, smooth, absolutely integral curve
of genus 1.  This is the curve we want to parameterize,
following the strategy presented in Sections~\ref{sec:mu3mu2} and~\ref{sec:selec}.
Since $r_s=r_t=2$ in this case, we look for  a map $B\rightarrow \PP^1$ of degree $3$ with two fully branched points and two
simply branched points. Such a map has two  totally ramified points.
They may be either $K$-rational or conjugated   over $K$. We will  assume that they
are $K$-rational. We call them  $P_0$ and $P_\infty$.
The two divisors $3P_0$ and $3P_\infty$ are  linearly equivalent because they both  are fibers of
the same  degree three map to $\PP^1$.
So  the difference $P_\infty -P_0$ has order $3$ in the Jacobian of $B$.
Our starting point will thus be a genus 1 curve $B/K$ and two points
$P_0$, $P_\infty$ in $B(K)$ such that $P_\infty-P_0$ has order $3$ in the Jacobian.

%\subsection{Constructing the parameterization}

Let $z\in K(B)$ be a function with divisor $3(P_0-P_\infty)$.
There is a unique hyperelliptic involution $\sigma : B\rightarrow B$ sending
$P_0$ onto $P_\infty$. It is defined over $K$.
There exists a scalar $a_{0,0}\in K^*$ such that
$\sigma(z)\times z=a_{0,0}$.
Let $x$ be a degree $2$ function, invariant by
$\sigma$, with polar divisor $(x)_\infty=P_0+P_\infty$.
Associated to the inclusion $K(x)\subset K(x,z)$ there is
a map $B\rightarrow \PP^1$ of degree $2$.
The sum $z+\sigma(z)$ belongs to $K(x)$. As a function on
$\PP^1$ it has a single pole of multiplicity $3$ at $x=\infty$.
So  $z+a_{0,0}/z$ is a polynomial of degree $3$ in $x$. Multiplying
$z$ by  a scalar, and adding a scalar to $x$, we may assume
that
\begin{equation}\label{eq:zxa}
z+\frac{a_{0,0}}{z}=x^3+a_{1,1}x+a_{0,1}.
\end{equation}
The image of $x\times z : B\rightarrow \PP^1\times \PP^1$ has
equation
\begin{equation*}\label{eq:eqgene1}
Z_0Z_1\left(X_1^3+a_{1,1}X_1X_0^2+a_{0,1}X_0^3\right)=X_0^3\left(Z_1^2+a_{0,0}Z_0^2\right).
\end{equation*}
This is a  curve $B^\star\subset \PP^1\times \PP^1$
with  arithmetic genus 2.
Since $B$ has geometric genus 1, we deduce that  $B^\star$ has one ordinary
double point (with finite $x$ and $z$ coordinates). Let
$(x,z)=(j,k)$ be this singular point. We find
\[a_{0,0}=k^2, \,\, a_{1,1}=-3j^2, \,\, a_{0,1}=2k+2j^3.\]
The  plane affine model $B^\star$ has equation
\begin{equation}\label{eq:gen1aff}
z^2+{k^2}=z\left(x^3-3j^2x+2(k+j^3)\right).
\end{equation}
This is a degree $3$ equation in $x$ with twisted discriminant
$81(1-k/z)^{2}$ times
\begin{eqnarray*}
h(z)&=&z^2-(2k+4j^3)z+k^2.
\end{eqnarray*}
We can parameterize $B$ with cubic radicals. We first parameterize the conic  $C$ with equation
\begin{equation}\label{eq:conic}
v^2=h(z)
\end{equation}
using the rational point $(z,v)=(0,k)$.
Applying  Tartaglia-Cardan formulae to the cubic  Equation~(\ref{eq:gen1aff})
we deduce a parameterization of $B$ with one cubic radical.
In order to relate Equation~(\ref{eq:gen1aff}) to a Weierstrass model, we simply sort in $z$ and find the degree $2$ equation
in $z$,
\[z^2   -(x^3-3j^2x+2k+2j^3)z+k^2=0\]
with discriminant
\[(x^3-3j^2x+2k+2j^3)^2-4k^2=(x-j)^2(x+2j)(x^3  -3j^2x+4k+2j^3).\]
A Weierstrass model for $B$ is then $u^2=(x+2j)(x^3  -3j^2x+4k+2j^3)$.
Replacing $j$ by $\lambda j$ and $k$ by $\lambda^3k$ for some
non-zero $\lambda$ in $K$  we obtain an isomorphic curve. So we may assume that
$j\in \{0,1\}$ without loss of generality.
This  construction is not substantially different from
the ones given by  Farashahi \cite{fara} and Kammerer, Lercier,
Renault \cite{KLR}.
Starting from any genus 1 curve $B$ and two points $P_0$ and $P_\infty$
such that $P_\infty-P_0$ has order $3$ in the Jacobian, we can construct a model of $B$
as in Equation~(\ref{eq:gen1aff}) and a parameterization of $B$.

\subsection{Example}

Let us consider an elliptic curve given in Weierstrass form $Y^2 = X^3+a\,X
+b$, for example  the curve $Y^2 = X^3 + 3\,X - 11$ over $\RR$,
together with a 3-torsion point $(x_0,y_0) = (3,-5)$.
Define the scalars   $\alpha$ and $\beta$ by
\begin{displaymath}
  \alpha = -\,{\frac {3\,{{x_0}}^{2}+a}{{2\,y_0}}}\text{ and }\beta = -y_0 - \alpha\,x_0.
\end{displaymath}
The functions $x=\alpha/3 + (Y+y_0)/(X-x_0)$ and $z=Y+\alpha\,X+\beta$
have divisors with zeros and poles as prescribed. On our particular  curve, these
functions are
\begin{equation}\label{eq:1}
  x = {\frac {Y-5}{X-3}} + 1\text{ and } z = Y + 3\,X-4\,.
\end{equation}
The functions $x$ and $z$ are related by Equation~(\ref{eq:zxa}) where
\begin{displaymath}
  a_{0,0}=4\,y_0^2=100\,,\ %
  a_{1,1}=-4\,x_0=-12\,,\ %
  a_{0,1}=-4\,{\frac{4\,{a}^{3}+27\,{b}^{2}}{{27\,{y_0}}^{3}}}=4\,.
\end{displaymath}
So
\[z+\frac{100}{z} = x^3-12\,x+4\,.\]
The double point on the latter is $(x,z)=(j,k)$ with
\begin{displaymath}
  j = \frac { -2\alpha}{ 3} = -2\text{ and }k =
-2\,y_0 = 10\,.
\end{displaymath}
A parameterization of the conic $C$ given by Equation~(\ref{eq:conic}) that reaches the point $(z,v)=(0,k)$ at
$t=\infty$ is
\begin{displaymath}
  z = 2\,{\frac {kt-k-2\,{j}^{3}}{{t}^{2}-1}} = 4\,{\frac {5\,t+3}{{t}^{2}-1}}\, ,\ %
  v = k-tz=\frac{(2k+4j^3)t-kt^2-k}{t^2-1}\, .
\end{displaymath}
and using Tartaglia-Cardan formulae  we find
 $x={\rho}/{3}+\, 3j^2/{\rho}$ with
\begin{displaymath}
  \rho =3j^2\times \sqrt [3] {\frac {2(t+1)}{\left(2\,{j}^{3}-kt+k \right)\left( 1-t \right)}}\,.
\end{displaymath}
It remains  to invert Eq.~\eqref{eq:1} in order to express  $X$ and $Y$ as
functions of $x$ and $z$, \textit{i.e.} as functions of the parameter $t$. For
$t=0$, we obtain in this way the point
\begin{displaymath}
  (X,Y)=(2\,(\sqrt [3]{3})^2+4\,\sqrt [3]{3}+3,\ -6\,(\sqrt [3]{3})^2-12\,\sqrt [3]{3}-17)\,.
\end{displaymath}

\section{Curves of genus 2}\label{sec:g2}

We look for parameterizations of
genus 2 curves. We will follow the strategy of Sections~\ref{sec:mu3mu2} and~\ref{sec:selec}.
We take  $r_s=4$ and $r_t=2$  this time.
Given a genus 2 curve $B$, we look for a degree three map $B\rightarrow \PP^1$
having $4$ simply branched points and $2$ totally branched points.
Such a map has two  totally ramified points.
We will  assume that they
are $K$-rational. We call them  $P_0$ and $P_\infty$.
The difference $P_\infty -P_0$ has order $3$ in the Jacobian of $B$.
Our starting point will thus be a genus 2  curve $B/K$ and two points
$P_0$, $P_\infty$ in $B(K)$ such that $P_\infty-P_0$ has order $3$ in the Jacobian.
The calculations will be slightly different depending on whether the set $\{P_0,P_\infty\}$
is stable under the action of the hyperelliptic involution of $B$ or not.
These two cases will be treated in Sections~\ref{sec:2d} and~\ref{sec:compl} respectively.
Section~\ref{sec:gene2} recalls simple facts about genus 2 curves.
Explicit calculations are detailed in Sections~\ref{sec:wg2} and~\ref{sec:exg2}.

\subsection{Generalities}\label{sec:gene2}

Let $K$ be a field of odd characteristic. Let
$\bar K$ be an algebraic closure of $K$.
Let $B/K$ be a projective, smooth, absolutely integral curve
of genus~2. Take  two
non-proportional holomorphic differential forms and let
$x$ be their quotient. This  is a function
on $B$ of degree $2$. Any degree $2$ function $y$ on $B$
belongs to the field $K(x)\subset K(B)$. Otherwise
the image of  $x\times y : B\rightarrow \PP^1\times \PP^1 $ would
be a curve birationally equivalent to $B$ with  arithmetic
genus $(2-1)\times (2-1)=1$.
A contradiction. So every degree two function on $B$
has the form $(ix+j)/(kx+l)$ with $i$, $j$, $k$ and $l$
in $K$. And $B$ has a unique hyperelliptic involution
$\sigma$. This is the non-trivial automorphism of the Galois
extension $K(x)\subset K(B)$. From Hurwitz genus formula, this
extension is ramified at exactly $6$ geometric points
$(P_i)_{1\le i\le 6}$
in  $B(\bar K)$.
If $\#K>5$ we can assume that the unique pole of
$x$ is not one of the $P_i$.
%Set $\alpha_i = x(P_i) \in \bar K$ for $1\le i\le 6$.
Set $F(x)=\prod_i(x-x(P_i))\in K[x]$.
According to Kummer theory, there exists a scalar $F_0 \in K^*$
such that $F_0F$
has a square root $y$ in $K(B)$. We set $f=F_0F$ and obtain
an  affine model for $B$ with  equation
\[y^2=f(x)\] and two $\bar K$-points $O$ and $\sigma(O)$
at infinity.
Every function $c$ in $K(B)$ can be written as
\[c=a(x)+yb(x)\]
with $a$ and $b$ in $K(x)$. If $P=(x_P,y_P)$ is a $\bar K$-point  on $B$ we denote
by $v_P$ the associated valuation of $\bar K(B)$. If $P$ is one of the
 $(P_i)_{1\le i\le 6}$  then
\begin{equation}\label{eq:v1}
v_P(c)=\min (2v_{x_P}(a),2v_{x_P}(b)+1),
\end{equation}
where $x_P=x(P)\in \bar K$ and $v_{x_P}$ is the valuation of $\bar K(x)$
at $x=x_P$.
If $P$ is a finite point which is not fixed by  $\sigma$ then
\begin{equation}\label{eq:v2}
\min(v_P(c), v_{\sigma(P)}(c))=
\min (v_{x_P}(a),v_{x_P}(b)).
\end{equation}
Finally
\begin{equation}\label{eq:v3}
\min(v_O(c), v_{\sigma(O)}(c))=\min (-\deg(a),-\deg(b)-3).
\end{equation}

Let $J$ be the Jacobian of $B$. A point $x$ in $J$ can be represented
by a divisor in the corresponding linear equivalence class.
We may fix a degree $2$ divisor $\Omega$ and associate to $x$
a degree $2$ effective divisor $D_x$ such that $D_x-\Omega$
belongs to the linear equivalence class associated with $x$.
This $D_x$ is generically unique. Indeed the only special
effective divisors  of degree $2$ are the fibers of the map
$B\rightarrow \PP^1$.
 We may also represent linear equivalence classes by divisors
of the form $P-Q$ where  $P$ and $Q$ are points on $B$.
There usually are two such representations as the map
\begin{displaymath}
    \xymatrixrowsep{0.1in}
    \xymatrixcolsep{0.3in}
\xymatrix{
    \xymatrixrowsep{0.1in}
    \xymatrixcolsep{0.3in}
 B^2 \ar@{->}[r]   & \Jac (B) \\
 (P,Q)  \ar@{|->}[r] & P-Q,
}
\end{displaymath}
is surjective and its restriction to the open set defined by
\[P\not =Q, P\not =\sigma(Q)\]
is finite {\'e}tale of degree $2$.

\subsection{A 2-dimensional family}\label{sec:2d}

Let $K$ be a field of  characteristic prime to $6$.
%Let $\bar K$ be an algebraic closure of $K$.
In this paragraph we study genus 2 curves $B/K$
satisfying the condition that
there exists a point $P$ in $B(K)$  such that the class of
$\sigma(P)-P$  has order $3$
in the Picard group.
In particular $P$ is not fixed by $\sigma$.
%Because this is a geometric condition we will assume that $K$ is algebraically closed.
We let $x$ and $y$ be  functions as in Section~\ref{sec:gene2}.
We can assume that  $x(P)=\infty$.
Let $z$  be a function with divisor $3(\sigma(P)-P)$. There exists a scalar $w\in K^*$ such that
$\sigma(z)\times z=w$.
We write
\[z=a(x)+yb(x)\]
with $a$ and $b$ in $K(x)$.
We deduce from Equations~(\ref{eq:v1}),  (\ref{eq:v2}),
 (\ref{eq:v3}), that $a$ and $b$ are polynomials
and $\deg(a)\le 3$ and $\deg(b)\le 0$.
From $z\sigma(z)=a^2-b^2f=w\in K^*$ we deduce that $\deg(b) = 0$
and $\deg(a) = 3$. We may divide  $z$  by a scalar in $K^*$ and
assume that $a$ is unitary. Replacing $x$ by $x+\beta$ for some
$\beta$ in $K$, we may even assume that $a(x)=x^3+kx+l$ with $k$ and
$l$ in $K$. Replacing  $y$ by
$b y$  we may assume that
 $b=1$ so \[z=y+x^3+kx+l.\]
An  affine  plane model  for $B$ has thus equation
\[z^2-2a(x)z+w=0\] that is
\begin{equation}\label{eq:cub1}
x^3+kx+l=\frac{z+wz^{-1}}{2}.
\end{equation}
This is a degree $3$ equation in $x$ with coefficients $s_1=0$, $s_2=k$, $s_3=(z+wz^{-1})/2-l$, and twisted discriminant
$81/4$ times
\[h(z)=z^2+w^2z^{-2} - 4l(z+wz^{-1})+2w+4l^2+\frac{16k^3}{27}.\]
We can parameterize
$B$ with cubic radicals. We first parameterize the elliptic curve $C$ with equation $v^2=h(z)$ with one cubic radical,
using e.g. Icart's method \cite{icart}. We deduce a parameterization
of $B$ applying Tartaglia-Cardan formulae to the cubic  Equation~(\ref{eq:cub1}). This introduces another cubic radical.
This is essentially the construction given by Kammerer, Lercier and Renault \cite{KLR}. Note that this
family of genus 2 curves has  dimension 2: when $K$ is algebraically closed  we may assume  that
$w=1$ without loss of generality.

\subsection{The complementary 3-dimensional family}\label{sec:compl}

We  still assume that $K$ %is algebraically closed and
has prime to $6$  characteristic.
We  consider a genus 2 curve $B$
and  two points $P_0$ and $P_\infty$ in $B(K)$
such that the difference $P_0-P_\infty$
has order $3$ in the Picard group.
This time we assume that $P_\infty\not=\sigma(P_0)$.
There exists
 a degree $2$ function $x$ having a zero at $P_0$ and a pole at $P_\infty$.
Let $z$ be a function with divisor $3(P_0-P_\infty)$.
The image of $x\times z : B\rightarrow \PP^1\times \PP^1$ has
equation
\begin{equation}\label{eq:eqgene}
\sum_{\substack{0\leqslant i\leqslant 3 \\0\leqslant j\leqslant 2}}a_{i,j}X_1^iX_0^{3-i}Z_1^jZ_0^{2-j}=0.
\end{equation}
The function $z$ takes the value  $\infty$ at  a single  point, and  $x$
has a  pole at this point. So if we set $Z_0=0$ in Equation~(\ref{eq:eqgene}) the form we  find  must be  proportional to
 $Z_1^2X_0^3$. We deduce that
\[a_{3,2}=a_{2,2}=a_{1,2}=0\] and
\[a_{0,2}\not =0.\]
The function  $z$ takes value  $0$ at a single  point, and  $x$
has a zero at this  point. So if we set  $Z_1=0$ in Equation~(\ref{eq:eqgene}) the form we find must be   proportional to
 $Z_0^2X_1^3$. We deduce that
\[a_{2,0}=a_{1,0}=a_{0,0}=0\] and
\[a_{3,0}\not =0.\]
Equation~(\ref{eq:eqgene}) now reads
\begin{equation*}\label{eq:casgp}
(a_{3,0}Z_0+a_{3,1}Z_1)Z_0X_1^3+(a_{1,1}X_0+a_{2,1}X_1)Z_0Z_1X_0X_1+
(a_{0,1}Z_0+a_{0,2}Z_1)Z_1X_0^3=0.
\end{equation*}
This is a  curve of arithmetic genus 2 in $\PP^1\times \PP^1$. It must be smooth because it  has geometric genus 2.
The corresponding plane affine model has equation

\begin{equation}\label{eq:casga}
(a_{3,0}+a_{3,1}z)x^3+(a_{1,1}+a_{2,1}x)zx+
(a_{0,1}+a_{0,2}z)z=0.
\end{equation}
This is a degree $3$ equation in $x$ with twisted discriminant
$z^{2}(a_{3,0}+a_{3,1}z)^{-4}$ times
\begin{eqnarray*}
h(z)&=&(9a_{0,2}a_{3,1})^2z^4+(12a_{0,2}a_{2,1}^3+162a_{3,0}a_{0,2}^2a_{3,1}-54a_{1,1}a_{2,1}a_{0,2}a_{3,1}\\&+&162a_{0,1}a_{3,1}^2a_{0,2})z^3
+(81a_{3,0}^2a_{0,2}^2+12a_{0,1}a_{2,1}^3-54a_{1,1}a_{2,1}a_{0,1}a_{3,1}\\&+&324a_{3,0}a_{0,1}a_{0,2}a_{3,1}-3a_{1,1}^2a_{2,1}^2-54a_{3,0}a_{1,1}a_{2,1}a_{0,2}+81a_{0,1}^2a_{3,1}^2\\&+&12a_{3,1}a_{1,1}^3)z^2+(12a_{1,1}^3a_{3,0}-
54a_{3,0}a_{1,1}a_{2,1}a_{0,1}+162a_{3,0}^2a_{0,1}a_{0,2}\\&+&162a_{3,0}a_{0,1}^2a_{3,1})z+(9a_{3,0}a_{0,1})^2.
\end{eqnarray*}
We can parameterize $B$ with cubic radicals. We first parameterize the elliptic curve with equation $v^2=h(z)$ with one cubic radical,
using Icart's method. We deduce a parameterization
of $B$ applying  Tartaglia-Cardan formulae to the cubic  Equation~(\ref{eq:casga}). This introduces another cubic radical.

In order to relate Equation~(\ref{eq:casga}) to a
hyperelliptic  model, we simply sort in $z$ and find the degree $2$ equation
in $z$,
\[a_{0,2}z^2+   (a_{3,1}x^3+a_{2,1}x^2+a_{1,1}x+a_{0,1}) z +  a_{3,0}x^3=0\]
with discriminant
\begin{equation}\label{eq:mx}
m(x)=(a_{3,1}x^3+a_{2,1}x^2+a_{1,1}x+a_{0,1})^2-4a_{0,2}a_{3,0}x^3.
\end{equation}
A hyperelliptic  model for $B$ is then \[y^2=m(x).\]
The construction will succeed for every genus 2 curve having a rational
3-torsion point in its
Jacobian that splits in the sense that it can be represented as a difference between two $K$-rational
points on $B$.

\subsection{Rational 3-torsion points in genus 2 Jacobians}\label{sec:wg2}

In this section we  start from a hyperelliptic  curve
\[\by^2=\bm(\bx),\] where  $\bm(\bx)$ is a degree
$6$ polynomial. We look for a  parameterization
 of it,  following  Sections~\ref{sec:2d}
or~\ref{sec:compl}. To this end  we need  a model
as in Equations~(\ref{eq:casga}) and~(\ref{eq:mx}). Such  a model is obtained  by
writing $\bm(\bx)$ as a difference $\bm_3(\bx)^2-\bm_2(\bx)^3$ where $\bm_3$ is a degree
$\le 3$ polynomial and $\bm_2$ is a degree
$\le 2$ polynomial with rational
roots.
We now are very close to investigations by Clebsh \cite{CLE} and Elkies
\cite{ELK}. Three-torsion points in the Jacobian of the curve
$\by^2=\bm(\bx)$ correspond
to expressions of $\bm$ as a difference between a square and a cube.
When the base field $K$ is finite, we may  first compute
 the Zeta function of the curve, deduce the cardinality of the
Picard group and obtain  elements of order~$3$ in it
by multiplying random elements
in the Picard group by the prime to three part of its order.
For a general base field $K$, we can look for solutions to
$\bm(\bx)=\bm_3(\bx)^2-\bm_2(\bx)^3$ by a direct Gr\"obner basis computation.
 Our experiments with the computer algebra softwares
\textsc{maple} or \textsc{magma} show that this
 approach is
efficient enough when $K$ is a finite field of
reasonable (say cryptographic) size.
When $K$ is the field $\QQ$ of rationals, this direct approach
becomes quite slow.

In this section we explain how to accelerate the computation
using invariant theory.
%We consider  the set of triples
%$(m(X,Y),m_2(X,Y),m_3(X,Y))$ where $m$, $m_2$, $m_3$ are homogeneous
%forms of respective degrees $6$, $2$, and $3$, such that $m=m_3^2-m_2^3$.
%This set   is invariant
%under linear change of
%variables.
Our
method takes as input, instead of $\bm(\bx)$, the
standard  homogeneous invariants for the action of $\GL_2$
evaluated at $\bm(\bX_1,\bX_0)$, the
degree $6$ projective form associated with $\bm(\bx)$.
Classical invariant theory results~\cite{Bolza,CLE} show
that the orbit under  $\GL_2$ of a
 degree 6 non-singular form  $\bm(\bX_1,\bX_0)$ is characterized  by
5 homogeneous  invariants $I_2$, $I_4$, $I_6$, $I_{10}$, $I_{15}$,
of respective degrees 2,
4, 6,  10, and 15.
There is a degree 30 algebraic relation  between the  $I_i$  (see~\cite{Igusa60}).

The action of $\GL_2$ on pairs $(\bm_2(\bX_1,\bX_0),\bm_3(\bX_1,\bX_0))$ consisting of
a quadric   and a cubic gives rise to well known invariants also:
$\iota_2$ (the discriminant of
$m_2$), $\iota_4$ (the discriminant of $m_6$) and 3 joint  invariants $\iota_3$,
$\iota_5$ and $\iota_7$, of respective degrees 2,  4, 3, 5 and 7. There is
 a degree 14 algebraic relation between the $\iota_i$~\cite[p.187-189]{Salmon1900}.
Since  the map $(\bm_2,\bm_3)\mapsto \bm=\bm_3^2-\bm_2^3$ is
$\GL_2$-equivariant we can describe its fibers in terms of the
invariants on each side.
We easily obtain the $I_i$'s as functions of the $\iota_i$'s,
\begin{eqnarray}\label{eq:grosys}
  2^2\,{I_2} &=& 120\,{\iota_5}+4\,{\iota_4}-12\,{\iota_3}\,{\iota_2}+3\,{{     \iota_2}}^{3}\,,\nonumber \\
  2^{7}\,{I_4} &=& 2640\,{{\iota_5}}^{2}+96\,{\iota_5}\,{\iota_4}-768\,{
    \iota_5}\,{\iota_3}\,{\iota_2}+240\,{\iota_5}\,{{\iota_2}}^{3}-24\,{\iota_4}\,{
    \iota_3}\,{\iota_2}+8\,{\iota_4}\,{{\iota_2}}^{3}\nonumber \\&&-8\,{{\iota_3}}^{3}+
  48\,{{
      \iota_3}}^{2}{{\iota_2}}^{2}-24\,{\iota_3}\,{{\iota_2}}^{4}+3\,{{\iota_2}}^
  {6}\,,\nonumber \\
  2^{10}\,{I_6} &=& -5120\,{{\iota_5}}^{3}-192\,{{\iota_5}}^{2}{\iota_4}-2304
  \,{{\iota_5}}^{2}{\iota_3}\,{\iota_2}+3504\,{{\iota_5}}^{2}{{\iota_2}}^{3}-
  96\,{\iota_5}\,{\iota_4}\,{\iota_3}\,{\iota_2}\nonumber \\&&+240\,{\iota_5}\,{\iota_4}\,{{
      \iota_2}}^{3}
  -288\,{\iota_5}\,{{\iota_3}}^{3}+1008\,{\iota_5}\,{{\iota_3}}^
  {2}{{\iota_2}}^{2}-768\,{\iota_5}\,{\iota_3}\,{{\iota_2}}^{4}+120\,{\iota_5
  }\,{{\iota_2}}^{6}\nonumber \\&&+4\,{{\iota_4}}^{2}{{\iota_2}}^{3}+24\,{\iota_4}\,{{
      \iota_3}}^{2}{{\iota_2}}^{2}
  -24\,{\iota_4}\,{\iota_3}\,{{\iota_2}}^{4}+4\,{
    \iota_4}\,{{\iota_2}}^{6}+36\,{{\iota_3}}^{4}{\iota_2}\\&&-72\,{{\iota_3}}^{3}{{
      \iota_2}}^{3}+48\,{{\iota_3}}^{2}{{\iota_2}}^{5}-12\,{\iota_3}\,{{\iota_2}}
  ^{7}+{{\iota_2}}^{9}\,,\nonumber \\
  2^{12}\,{I_{10}} &=& 46656\,{{\iota_5}}^{5}+3456\,{{\iota_5}}^{4}{\iota_4}-
  3888\,{{\iota_5}}^{4}{\iota_3}\,{\iota_2}+729\,{{\iota_5}}^{4}{{\iota_2}}^{
    3}+64\,{{\iota_5}}^{3}{{\iota_4}}^{2}\nonumber \\&&-144\,{{\iota_5}}^{3}{\iota_4}\,{
    \iota_3}\,{\iota_2}
  +27\,{{\iota_5}}^{3}{\iota_4}\,{{\iota_2}}^{3}+128\,{{\iota_5
    }}^{3}{{\iota_3}}^{3}-27\,{{\iota_5}}^{3}{{\iota_3}}^{2}{{\iota_2}}^{2}\,.\nonumber
\end{eqnarray}

Given the $I_i$'s evaluated at $\bm(\bX_1,\bX_0)$, the generic change of variable $\lambda
= {\iota_2}^3$ and $\mu = \iota_2\times \iota_3$ turns these equations into a
system of 4 equations of total degrees 1, 3, 4 and 6 in the 4 variables
$\lambda$, $\mu$, $\iota_4$ and $\iota_5$\,.  A Gr\"obner basis
can be easily computed for the lexicographic order (note that the first
equation is  linear). This yields a degree 40 polynomial in $\lambda$.
If none of the roots of this polynomial are squares, we can abort the
calculation because  we need
 $\bm_2(\bx)$ to have rational roots in order
to parameterize the curve $\by^2 = \bm(\bx)$.

Considering Equation~(\ref{eq:mx}) of Section~\ref{sec:compl}
it is natural to look for a form $m$ in the
$\GL_2$-orbite of $\bm$  such that
$m=m_3^2-m_2^3$ for some
 $m_2(x) = e\,x$ and $m_3(x) = a\,x^3 + b\,x^2 + c\,x + d$, where
$e^3=4a_{0,2}a_{3,0}$, $a=a_{3,1}$, $b=a_{2,1}$, $c=a_{1,1}$, $d=a_{0,1}$.
The invariants of $(m_2,m_3)$ are
\begin{eqnarray}\label{eq:2}
  &&\iota_2 = {e}^{2}\,,\ \iota_3 = -e ( 9\,ad-bc )\,,\ \ \ \ \ \ \ \ \ \iota_5 = -{e}^{3}ad\,,\ \iota_7 = {e}^{3} ( a{c}^{3}- {b}^{3}d )\,,\\
  &&\hspace*{2cm}\iota_4 = -27\,{a}^{2}{d}^{2}+18\,abcd-4\,a{c}^{3}-4\,{b}^{3}d+{b}^{2}{c}^{2}\,.\nonumber
\end{eqnarray}
So for each candidate $(\iota_2, \iota_3, \iota_4, \iota_5)$ issued from
Equations~(\ref{eq:grosys}), we
 invert Eq.~\eqref{eq:2}. A Groebner basis
for the lexicographic order $d$, $c$, $b$, $a$, $e$ yields generically a
1~dimensional system the  last two equations of which are
\begin{eqnarray*}
  0 &=& e^2 - \iota_2\,,\\
  0 &=& {{\iota_2}}^{3}{\iota_5}\,{b}^{6}-{\iota_2}\, ( {{\iota_2}}^{3}{\iota_4}-{{\iota_2}}^{2}{{\iota_3}}^{2}+36\,{\iota_2}\,{\iota_3}\,{\iota_5}-216\,{{\iota_5}}^{2} )\,e\,a\, {b}^{3}-4\, ( {\iota_2}\,{\iota_3}-9\,{\iota_5} ) ^{3}\,{a}^{2}\,.
\end{eqnarray*}
We  keep solutions $m_2(x)$ and $m_3(x)$ that
yield a polynomial $m(x) = m_3(x)^2-m_2(x)^3$ which is $\GL_2$-equivalent
 to
$\bm(\bx)$ over the base field (see \cite{LRS} for efficient algorithms). Applying
the isomorphism to $m_2(x)$ and $m_3(x)$ gives  $\bm_2(\bx)$ and $\bm_3(\bx)$.

\subsection{An example}\label{sec:exg2}

Let $K$ be a field with $83$ elements.
We start from
 the genus 2 curve with affine equation $\by^2=\bm(\bx)$ with
$\bm(\bx)=\bx^6+39\bx^5+64\bx^4+7\bx^3+\bx^2+19\bx+36$.
In order to find $\bm_3(\bx)$ and $\bm_2(\bx)$ such that $\bm(\bx) = \bm_3(\bx)^2 - \bm_2(\bx)^3$, we first
compute the invariants of the degree six form $\bm$
\begin{displaymath}
  (I_2, I_4, I_6, I_{10}) = (23, 9, 38, 53, 59)\,.
\end{displaymath}
A Groebner basis for the relations between
$\lambda$, $\mu$ and $\iota_4$ is
\begin{eqnarray*}
  \iota_4 &=&
  27\,{\lambda}^{39}+58\,{\lambda}^{38}+3\,{\lambda}^{37}+18\,{\lambda}^{36}+42\,{\lambda}^{35}+26\,{\lambda}^{34}+52\,{\lambda}^{33}+60\,{\lambda}^{32}\\&&+78\,{\lambda}^{31}+17\,{\lambda}^{30}+
  50\,{\lambda}^{29}+12\,{\lambda}^{28}+75\,{\lambda}^{27}+20\,{\lambda}^{26}+75\,{\lambda}^{25}+38\,{\lambda}^{24}\\&&+19\,{\lambda}^{23}+21\,{\lambda}^{22}+35\,{\lambda}^{21}+31\,{\lambda}^{20}+27\,{\lambda}^{19}+49\,{\lambda}^{18}+44\,{\lambda}^{17}+30\,{\lambda}^{16}\\&&+38\,{\lambda}^{15}+55\,{\lambda}^{14}+59\,{\lambda}^{13}+6\,{\lambda}^{12}+2\,{\lambda}^{11} +36\,{\lambda}^{10}+18\,{\lambda}^{9}+2\,{\lambda}^{8}+41\,{\lambda}^{7}\\&&+62\,{\lambda}^{6}+3\,{\lambda}^{5}+49\,{\lambda}^{4}+{\lambda}^{3}+33\,{\lambda}^{2}+36\,\lambda+69\,,\\
\mu &=&
62\,\lambda^{40}+46\,\lambda^{39}+11\,\lambda^{38}+33\,\lambda^{37}+75\,\lambda^{36}+19\,\lambda^{35}+53\,\lambda^{34}+10\,\lambda^{33}\\&&
+48\,\lambda^{32}+47\,\lambda^{31} +77\,\lambda^{30}+14\,\lambda^{29}+49\,\lambda^{28}+47\,\lambda^{27}+38\,\lambda^{26}+19\,\lambda^{25}\\&&+25\,\lambda^{24}+44\,\lambda^{23}+68\,\lambda^{22}+15\,\lambda^{21}+36\,\lambda^{20}+9\,\lambda^{19}+73\,\lambda^{18}+13\,\lambda^{17}\\&&+64\,\lambda^{16}+5\,\lambda^{15}+67\,\lambda^{14}+82\,\lambda^{13}+69\,\lambda^{12}+9\,\lambda^{11}+69\,\lambda^{10}+35\,\lambda^{9}\\&&+57\,\lambda^{8}+57\,\lambda^{7}+7\,\lambda^{6}+11\,\lambda^{5}+37\,\lambda^{4}+78\,\lambda^{3}+10\,\lambda^{2}+73\,\lambda\,,\\
0&=&{\lambda}^{40}+48\,{\lambda}^{39}+67\,{\lambda}^{38}+35\,{\lambda}^{37}+50\,{\lambda}^{36}+23\,{\lambda}^{
  35}+4\,{\lambda}^{34}+12\,{\lambda}^{33}\\&&+37\,{\lambda}^{32}+49\,{\lambda}^{31}+40\,{\lambda}^{30}+71\,{\lambda}^{29}+60\,{\lambda}^{28}+79\,{\lambda}^{27}+19\,{\lambda}^{26}+81\,{\lambda}^{25}\\&&+82\,{\lambda}
^{24}+26\,{\lambda}^{23}+9\,{\lambda}^{22}+19\,{\lambda}^{21}+82\,{\lambda}^{20}+40\,{\lambda}^{19}+50\,{\lambda}^{18}+67\,{\lambda}^{17}\\&&+80\,{\lambda}^{16}+29\,{\lambda}^{15}+73\,{\lambda}^{14}+38\,{\lambda}^{13}+81\,{\lambda}^{12}+73\,{\lambda}^{11}+5\,{\lambda}^{10}+14\,{\lambda}^{9}\\&&+82\,{\lambda}^{8}+46\,{\lambda}^{7}+62\,{\lambda}^{6}+32\,{\lambda}^{5}+17\,{\lambda}^{4}+74\,{\lambda}^{3}+15\,{\lambda}^{2
}+30\,\lambda+43\,.
\end{eqnarray*}
Here, we only have two rational candidates for $(\lambda,\mu, \iota_4)$, the first
one  gives
\begin{displaymath}
(\iota_2, \iota_3, \iota_4, \iota_5) = (17, 51, 35, 55)\,.
\end{displaymath}
Now, inverting Eq.~\eqref{eq:2} yields 4 possibilities, all parameterized by $a$:
\begin{enumerate}
\item $\{\ {d}+74\,{{c}}^{3}=0,{c}\,{b}+45=0,{c}\,{a}+63\,{{b}}^{2}=0,{{b}}^{3}+23\,{a}=0,{e}+73=0\ \}$\,,
\item or $\{\ {d}+65\,{{c}}^{3}=0,{c}\,{b}+45=0,{c}\,{a}+73\,{{b}}^{2}=0,{{b}}^{3}+46\,{a}=0,{e}+73=0\ \}$\,,
\item or $\{\ {d}+18\,{{c}}^{3}=0,{c}\,{b}+38=0,{c}\,{a}+73\,{{b}}^{2}=0,{{b}}^{3}+37\,{a}=0,{e}+10=0\ \}$\,,
\item or $\{\ {d}+9\,{{c}}^{3}=0,{c}\,{b}+38=0,{c}\,{a}+63\,{{b}}^{2}=0,{{b}}^{3}+60\,{a}=0,{e}+10=0\ \}$\,.
\end{enumerate}
A solution to the first set of equations is, for $a=1$,
\begin{displaymath}
  m_3(x) = {x}^{3}+46\,{x}^{2}+73\,x+47\ \text{ and }\ m_2(x) = 10\,x,
\end{displaymath}
and the polynomial
\[m(x)= m_3(x)^2-m_2(x)^3\]
is $\GL_2$-equivalent to $\bm(x)$. Indeed
\[m(\frac{76\bx+70}{36\bx+43})\times (36\bx+43)^6=\bm(\bx).\]
So we set
\[\bm_3(\bx)=m_3(\frac{76\bx+70}{36\bx+43})\times (36\bx+43)^3=15\bx^3 + 30\bx^2 + 46\bx + 7\]
and
\[\bm_2(\bx)=m_2(\frac{76\bx+70}{36\bx+43})\times (36\bx+43)^2=53\bx^2 + 29\bx + 54\]
and we check that $\bm=\bm_3^2-\bm_2^3$.

\subsection{Parameterization}

The curve  with equation  $\by^2=\bm(\bx)$
over the field with $83$ elements is isomorphic
to the curve with equation
\[y^2 = (a{x}^{3}+b\,{x}^{2}+c\,x+d)^2 - (e\,x)^3 = ({x}^{3}+46\,{x}^{2}+73\,x+47)^2 - (10\,x)^3\]
through the change of variables
\begin{equation}\label{eq:xX}
x=\frac{76\bx+70}{36\bx+43}, \text{ and }  y = \frac{\by}{(36\bx+43)^3}.
\end{equation}
With the notation in Section~\ref{sec:compl} we have
$a=a_{3,1}=1$, $b=a_{2,1}=46$, $c=a_{1,1}=73$,
$d=a_{0,1}=47$, $e=10$, $a_{0,2}=-1/2$, $a_{3,0}=-e^3/2$.
Let $P_0$ be the point with coordinates $x=0$ and $y=-47$. Let
$P_\infty$ be the point where $x$ has a pole and $y/x^3=1$.
The functions $x$ has a zero at $P_0$ and a pole
at $P_\infty$. The function  $z = y + a\,x^3 + b\,x^2 + c\,x + d$ has
 divisor $3(P_0-P_\infty)$. These two functions are related by the equation
\begin{equation}\label{eq:3}
  ( -{e}^{3}/2+az ) {x}^{3}+ ( bx+c ) zx + ( d-z/2 ) z = 0,
\end{equation}
that is  $( z+81 ) {x}^{3}+ ( 46\,x+73 ) zx+ ( 47+
41\,z ) z = 0$\,.
%
%With the notation above we have $a=a_{3,1}=1$,
%$b=a_{2,1}=46$, $c=a_{1,1}=73$, $d=a_{0,1}=47$, $a_{0,2}=-1/2=41$,
%$a_{3,0}=-e^3/2=81$.
The resolvent elliptic curve has equation $v^2 = h(z)$ with
\begin{displaymath}
  h(z) = 41\,{z}^{4}+15\,{z}^{3}+38\,{z}^{2}+46\,z+7\,.
\end{displaymath}
It
 is birationally isomorphic to the Weierstrass
curve with equation
\begin{math}
  Y^2 = X^3 + 37\,X + 60\,,
\end{math}
whose  Icart's parameterization in $t$  is
\begin{displaymath}
  X = \kappa/6+\,{t}^{2}/3\,,\ Y = ({t}^{3}+{t\,\kappa}+28/{t})/6
\end{displaymath}
where
\begin{displaymath}
  \kappa = \sqrt[3]{81\,{t}^{6}+79\,{t}^{2}+71+\frac{56}{{t}^{2}}}\,.
\end{displaymath}
After a birational change of variable, we obtain
\begin{eqnarray*}
z &=& {\frac {10\,Y+16\,X+72}{74\,{X}^{2}+79\,X+49}},\\
v &=& {\frac { ( 47\,{X}^{2}+8\,X+64 ) Y+51\,{X}^{4}+5\,{X}^{3}+
20\,{X}^{2}+20\,X+18}{81\,{X}^{4}+72\,{X}^{3}+47\,{X}^{2}+23\,X+77}}\,.
\end{eqnarray*}
We then apply Tartaglia-Cardan formulae to
Eq.~\eqref{eq:3} in order to obtain $x$ and $y=z-m_3(x)$ as functions of $t$.
Inverting the change of variables in Equation~(\ref{eq:xX})
gives a point   $(\bx,\by)$ on the initial curve.
%At $t=1$, for instance, we find so the point $(13, 52)$.

\subsection{The density of target curves}

We prove that the construction in Section~\ref{sec:compl} provides
a parameterization for a fixed positive proportion of genus 2 curves
over $\FF_q$ when $q$ is prime to $6$ and large enough.
We call $\cS$ the set of  non-degenerate sextic binary forms
with coefficients in $\FF_q$.
Scalar multiplication  \[(\lambda,m(X_1,X_0))\mapsto \lambda m(X_1,X_0)\]
defines an action of the  multiplicative group $\Fqs$
on $\cS$.
The linear group $\GL_2(\Fq)$ also acts on $\cS$.
Call $G$  the subgroup of $\GL_2(\Fq)\times \Fqs$
consisting of pairs $( \gamma, \lambda)$ where $\lambda$
is a square.
To every non-degenerate sextic binary form $m(X_1,X_0)$
with coefficients in $\FF_q$ we associate the $\FF_q$-isomorphism class
of the curve
with equation $y^2=m(x,1)$. This defines a surjective
map $\nu$ from $\cS$
onto the set $\cI$ of  $\FF_q$-isomorphism
classes of genus 2 curves over $\Fq$.
The fibers of $\nu$ are
the orbites for the action of  $G$ on $\cS$.
When $q$ tends to infinity, the proportion of forms in $\cS$
with non-trivial stabilizer in $G$ tends to zero.
So  it is equivalent to count
isomorphism classes of curves in $\cI$ or to count  forms in $\cS$.

We call $\cP$ the set of pairs $(m_2,m_3)$ consisting of
a split quadratic form \[m_2(X_1,X_0)=(aX_1-bX_0)(cX_1-dX_0)\]
and  a cubic form $m_3$, such that $m_3^2-m_2^3$ is a
non-degenerate sextic form. The cardinality of $\cP$
is $q^7\times (1/2+o(1))$ when $q$ tends to infinity.
Let $\chi : \cP\rightarrow \cS$ be the map that sends
$(m_2,m_3)$ onto $m_3^2-m_2^3$. According to
work by Clebsh \cite{CLE}
and Elkies \cite[Theorem 3]{ELK}, fibers of $\chi$ have no more
than $240$ elements. So the image of $\chi$ has cardinality
at least $q^7\times (1/480 +o(1))$ and density at least $1/480+o(1)$.

\begin{theorem}
Let  $q$ be a   prime power
that is  prime to $6$.
The  proportion of all  genus 2 curves over the field with $q$ elements
that can be parameterized by $3$-radicals is at least $1/480+\epsilon(q)$
where $\epsilon$ tends to zero when $q$ tends to infinity.
\end{theorem}

\section{Other families of covers}\label{sec:other}

In Sections~\ref{sec:g1} and~\ref{sec:g2} we have studied two families of $\mu_3\rtimes \mu_2$ covers
corresponding to $(r_s,r_t)=(2,2)$ and  $(r_s,r_t)=(4,2)$ respectively.
In this section we quickly review a few other possibilities. We also present an interesting family  of
$\mu_5\rtimes \mu_2$ covers.

\subsection{The case $(r_s,r_t)=(4,1)$}

Both $B$ and $C$ have genus 1. The map $B\rightarrow E$ is any degree three map having a triple pole. If
$B$ is given by a  Weierstrass model, then for every scalar $t$, the function $y+tx$ will do. So we obtain a one
parameter family of parameterization of $B$ by elliptic curves $C_t$.
The resolvents
$C_t$ form a  non-isotrivial family. However,
we observed that the 3-torsion group scheme  $C_t[3]$ is isomorphic to $B[3]$ for every value of $t$.

\subsection{The case $(r_s,r_t)=(6,1)$}\label{sec:61}

Both $B$ and $C$ have genus 2. The map $B\rightarrow E$ is any degree three map having a triple pole. There
is one such map for every non-Weierstrass point $P$ on $B$.
We obtain a one
parameter family of parameterization of $B$ by genus 2 curves $C_P$.
The resolvents $C_P$ form a  non-isotrivial family. However,
we observed that the 3-torsion group scheme  $J_{C_P}[3]$ is isomorphic
to $J_B[3]$ for every  $P\in B$.

\subsection{The case $(r_s,r_t)=(8,1)$}

Both $B$ and $C$ have genus 3. The map $B\rightarrow E$ is a degree three map having a triple pole $P$. This
pole is a rational Weierstrass point. The curve $C$ is hyperelliptic.
For every genus 3 curve $B$ having a rational  Weierstrass point, we thus obtain a  parameterization of $B$ by a hyperelliptic  curve of genus 3.
Conversely, for every
 hyperelliptic curve of  genus 3 which we can parameterize, we obtain a parameterization for a 1-dimensional family of  non-hyperelliptic
genus 3 curves.

\subsection{Curves with a $\mu_5\rtimes \mu_2$ action}\label{sec:mu5mu2}

This time we assume that the characteristic of $K$ is prime to $10$.
Let $\zeta_5\subset \bar K$ be a primitive $5$-th root of unity.
We denote by $\mu_5\rtimes \mu_2$ the subgroup scheme of $\Sym (\mu_5)$ generated by $x\mapsto x^{-1}$
and $x\mapsto \zeta_5x$.
Let $A$ be a projective, absolutely integral,
smooth curve  over $K$. We assume that
 $\Aut(C\otimes_K\bar K)$ contains  the finite {\'e}tale $K$-group scheme
$\mu_5\rtimes \mu_2$.
We set  $B=A/\mu_2$, and $C=A/\mu_5$.
If   $C$ admits a  parameterization  by $S$-radicals as
in Equation~(\ref{eq:para}), and if
 the normalization $D'$ of the fiber product of $A$ and $D$ above $C$ is  absolutely integral, then
we can construct an $S\cup \{5\}$-parameterization of $B$ just as in Section~\ref{sec:mu3mu2}.
We assume that
$E=A/(\mu_5\rtimes \mu_2)$ has genus 0.
Let $r_d$ be the number of  branched points with ramification type $2,2,1$.  Let $r_t$ be the the number  of  branched points with ramification type $5$.
According to the Hurwitz Genus Formula \cite[III.4.12, III.5.1]{Stich} the genus of $B$ is
\[g_B=r_d+2r_t-4.\]
Every  branched point of type $2,2,1$ of the cover $B\rightarrow E$ gives rise to a branched point
of type $2,2,2,2,2$ of the cover $A\rightarrow E$ and to a simple branched point of $C\rightarrow E$.
And every totally  branched point of the cover $B\rightarrow E$ gives rise to a branched point
of type $5,5$ of the cover $A\rightarrow E$ and to a non-branched point of $C\rightarrow E$. So
\[g_A=\frac{5r_d}{2}+4r_t-9, \text{ \, and \,\,\,}g_C=\frac{r_d}{2}-1.\]
We still call \[m=r_d+r_t-3\]  the {\it modular dimension}.
The { genericity condition}  is
\[2r_d+5r_t\le 12-2\epsilon(r_d+2r_t-4),\]
where $\epsilon(0)=3$, $\epsilon(1)=1$, and $\epsilon(n)=0$ for $n\ge 2$.

An interesting case is when $r_d=6$ and $r_t=0$. Then both $B$ and $C$ have genus~2. The map  $B\rightarrow E$ is a $\mu_5\rtimes \mu_2$-cover.
 The cover $A\rightarrow C$ is  unramified. It is a quotient
by $\mu_5$. Associated to it, there is a $C_5$ inside $J_C$.
So we are just dealing with a genus~2 curve $C$ having a 5-torsion
point in its Jacobian. We provide explicit equations for this situation
in Section~\ref{sec:g25}.

\section{Genus 2 curves with a 5-torsion divisor}\label{sec:g25}

We assume that  $K$ has characteristic prime to $10$. Let
$C$ be a genus 2 curve having a $K$-rational point of order~$5$ in its
Jacobian.
We assume that this point is the class of $P_\infty-P_0$ where
$P_\infty$ and $P_0$ are two $K$-rational points on $C$.
We give explicit equations for $C$, $P_0$ and $P_\infty$ depending
on rational parameters.
In Sections~\ref{sec:cas1}, \ref{sec:cas2}, and \ref{sec:cas3}, we
distinguish  three cases depending on the action of the hyperelliptic involution
$\sigma$ on  $P_0$ and $P_\infty$. We note that these two
points cannot be both Weierstrass points.
We finally give in Section~\ref{sec:exg25} an example of how to
combine this  construction and the previous ones
 in order to parameterize more genus 2 curves.

\subsection{A first special case}\label{sec:cas1}

We   first assume that $P_0$ is a  Weierstrass point. So
$P_\infty$ is not. Let $x$ be a degree $2$ function having a
pole at $P_\infty$ and a zero at $P_0$.
Let $y$ be a function as in Section~(\ref{sec:gene2}).
We have $y^2=f(x)$ for some degree $6$ polynomial in $K[x]$.
Let $z\in K(C)$ be a function with divisor $5(P_0-P_\infty)$.
We write \[z=a(x)+yb(x)\] with $a(x)$ and $b(x)$ in $K(x)$.
We deduce from Equations~(\ref{eq:v1}),  (\ref{eq:v2}),
 (\ref{eq:v3}), that $a$ and $b$ are polynomials
and $\deg(a)\le 5$ and $\deg(b)\le 2$.
Since $z$ has a pole of order $5$ at $P_\infty$ and has valuation
$0$ at $\sigma(P_\infty)$ we actually know that
$\deg(a) =  5$ and $\deg(b)= 2$. Also $b$ is divisible by $x$ exactly twice,
and $a$  is divisible by $x$ at least thrice. Multiplying $z$ by a scalar
we may ensure that $a$ is unitary. Multiplying $y$ by a scalar
we may ensure that $b=x^2$.  And  $a(x)=x^3(x^2+kx+l)$ for some  $k$ and some
$l$ in $K$. There exists a scalar $w\in K^*$ such that
\[z\times \sigma(z)=wx^5=x^4(x^2(x^2+kx+l)^2-f(x)).\]
So $f(x)=x^2(x^2+kx+l)^2-wx$. The curve $C$ has affine equation
\[y^2=x^2(x^2+kx+l)^2-wx,\]
$P_\infty$ is one  of the two points at infinity,
 and $P_0$ is the point $(0,0)$.
This is essentially the model  given by Boxall, Grant and Lepr{\'e}vost \cite{BGL}.

\subsection{Another special case}\label{sec:cas2}

We assume now that  $\sigma (P_0)=P_\infty$. Let $x$ be a degree two
function having  poles at $P_0$ and $P_\infty$.
Let $y$ and $f(x)$ be as in Section~\ref{sec:gene2}.
Let $z$ be a function
with divisor $5(P_0-P_\infty)$. We write $z=a(x)+yb(x)$ where
$a$ and $b$ are polynomials in $x$ with degrees $5$ and $2$.
 Multiplying $z$ by a constant in $K$ we may assume that $a$
is unitary. Multiplying $y$ by a constant in $K$ we may assume that $b$
is unitary.
Adding a constant to $x$ we may assume that
\[b(x)=x^2-k\] for some $k\in K$.
There is a scalar $w\in K^*$ such that
\[z\times \sigma(z)=w=a^2-fb^2.\]
So $w$ is a square in the algebra $K[x]/b(x)$. This leaves two
possibilities.
Either $w=W^2$ for some $W\in K^*$ and $a(x)=W\bmod b(x)$, or
$w=W^2k$  for some $W\in K^*$ and $a(x)=Wx\bmod b(x)$.
We study these two subcases successively.

\subsubsection{If $w=W^2$ and $a(x)=W\bmod b(x)$}
We check that \[a(x)=W\bmod b(x)^2\] indeed. Since $a$ is unitary,  there exists
a scalar $j\in K$ such that
$a=W+(x+j)b^2$.  We deduce expressions for $a$, $b$
and $f$ in the parameters $k$, $W$,  and $j$.
The actual dimension of the family is $2$ because we may multiply $x$
by a scalar.

\subsubsection{If $w=W^2k$  and $a(x)=Wx\bmod b(x)$}
In particular $k$ is not $0$. We check that
$a(x)=Wx+a_1(x)b(x)\bmod b(x)^2$ with $a_1(x)=-Wx/(2k)$. So there exist
a scalar  $j\in K$ such that
\[a=Wx -Wxb(x)/(2k)+(x+j)b(x)^2.\]  We deduce expressions for $a$, $b$
and $f$ in the parameters $k$, $W$,  and $j$.
The actual dimension of the family is $2$ again.

\subsection{Generic case}\label{sec:cas3}
We   assume
that none of $P_0$ and $P_\infty$
is a  Weierstrass point and $\sigma(P_0)\not =P_\infty$.
Let $x$ be a degree $2$ function having a zero at $P_0$ and a pole
at $P_\infty$. Let $y$ be a function as in Section~\ref{sec:gene2}.
We have $y^2=f(x)$ where $f\in K[x]$ is a degree $6$ polynomial.
Both $f(0)$ and the leading coefficient of $f$ are  squares in $K$.
Let $z\in K(C)$ be a function with divisor $5(P_0-P_\infty)$.
We write \[z=a(x)+yb(x)\] with $a(x)$ and $b(x)$ in $K(x)$.
We deduce from Equations~(\ref{eq:v1}),  (\ref{eq:v2}),
 (\ref{eq:v3}), that $a$ and $b$ are polynomials
and $\deg(a)\le 5$ and $\deg(b)\le 2$.
Since $z$ has a pole of order~$5$ at $P_\infty$ and has valuation
$0$ at $\sigma(P_\infty)$ we actually know that
$\deg(a) =  5$ and $\deg(b)= 2$.
Multiplying  $z$ by  a scalar, we may ensure that
$a$ is unitary.
Multiplying  $y$ by  a scalar, we may ensure that
$b$ is unitary.
Since $z$ has a zero of order~$5$ at $P_0$ and has valuation
$0$ at $\sigma(P_0)$ we know that
$a(0)\not = 0$ and $b(0)\not =0$.

The three polynomials $a(x)$, $b(x)$,
and $f(x)$ are related by the equation \[a^2-fb^2=wx^5\]
for some $w\in K^*$.
In particular, $wx$ is a square modulo $b(x)$.
We can easily deduce that
\[b(x)=x^2+(2k-wl^2)x+k^2\]
for some $k$ and $l$ in $K^*$. A square root of $wx$ modulo $b(x)$ is then
$(k+x)/l$. A square root of $wx^5$ modulo $b$ is then
\[a_0(x)=\frac{(k^2-3kwl^2+w^2l^4)x+(k-wl^2)k^2}{l}.\]
Using Hensel's lemma we deduce that $a$ is the
square root of $wx^5$ modulo $b^2$
of the form $a_0+a_1b$ with
\[a_1(x)=\frac{k^2-2wkl^2+2w^2l^4+x(wl^2+k)}{2wl^3}.\]

So there exists  $j\in K$ such that
$a=a_0+a_1b+a_2b^2$ with \[a_2(x)=x+j.\]
We deduce the expressions for $a$, $b$
and $f=(a^2-x^5)/b^2$ in the parameters $j$, $k$, $l$, $w$.
%Multiplying $x$ by a scalar we may assume that $k=1$.

\subsection{An example}\label{sec:exg25}

Let $K$ be a field with $83$ elements.
We set $w=1$, $j=2$, $k=3$, $l=14$ and find
$a(x)=x^5 + 37x^4 + 78x^3 + 18x^2 + 26x + 29$
and
$b(x)=x^2+59x+9$, and \[f(x)=x^6 + 39x^5 + 64x^4 + 7x^3 + x^2 + 19x + 36.\]
The curve $C$ with equation $y^2=f(x)$ has genus 2. Its Jacobian
has $3.5.7.71$ points over $K$. We set $z=a(x)+yb(x)$ and define
a cyclic unramified covering $A$ of $C$ by setting
$t^5=z$. We lift the action of the hyperelliptic involution
 $\sigma$ onto $A$ by setting
$\sigma(t)=x/t$. The function $u=t+x/t$ is invariant by $\sigma$.
The field $K(u,x)$ is the function field of the quotient curve
$B=A/\sigma$. A singular plane model for $B$ is given by the equation
\[u^5+78xu^3+5x^2u=2a(x)=2(x^5+37x^4+78x^3+18x^2+26x+29).\]
Note the Tchebychev polynomial on the left hand side.
The Jacobian of $B$ has $5.37^2$ points over $K$. In particular,
its 3-torsion is trivial. However we can parameterize the
curve $B$ using the parameterization of $C$ constructed
in Section~\ref{sec:exg2}. Note that $C$ appears in Section~\ref{sec:exg2}
under the name $B$.

\subsection{Composing parameterizations}

In Section~\ref{sec:exg25} we parameterize a genus 2 curve
(call it $B_2$) by
another genus 2 curve (call it $C_2$),
using a $\mu_5\rtimes \mu_2$ action on some curve $A_2$. In Section~\ref{sec:exg2}
we had constructed a parameterization of $C_2=B_1$ by a genus
one curve (call it $C_1$) using a $\mu_3\rtimes \mu_2$ action
on some curve $A_1$. This $C_1$ can be parameterized e.g.
using Icart's parameterization. Composing the three parameterizations
we obtain a parameterization of $B_2$ by $\PP^1$.

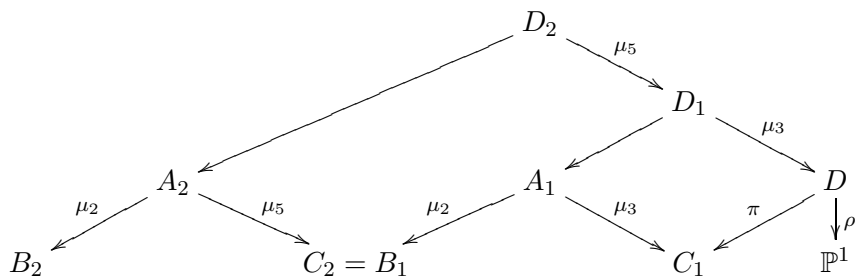
\begin{figure}\centering
  \begin{equation*}
    \xymatrixcolsep{0.5in}\xymatrixrowsep{0.2in}\xymatrix{&&&D_2\ar@{->}[ddll]\ar@{->}[dr]^{\mu_5}&&\\
      &&&&D_1 \ar@{->}[dl]\ar@{->}[dr]^{\mu_3} &   \\
      &A_2\ar@{->}[dl]_{\mu_2}\ar@{->}[dr]^{\mu_5}&&A_1\ar@{->}[dl]_{\mu_2}\ar@{->}[dr]^{\mu_3}&&D\ar@{->}[dl]_{\pi}\ar@{->}[d]^{\rho}\\
      B_2&&C_2=B_1&&C_1&\PP^1
    }\end{equation*}
\caption{Composing parameterizations}\label{fig:passage3}
\end{figure}

This situation is represented on Figure~\ref{fig:passage3}.
The curve $D_1$
is the normalization of the fiber product of $D$ and $A_1$ over $C_1$. The curve
$D_2$ is the normalization of the fiber product of $D_1$ and $A_2$ over $C_2$.
 We can prove that $D_1$ and
$D_2$ are  absolutely irreducible by observing that all down left arrows
have degree a power of two, while all down right arrows are Galois
of  odd degree.
The interest of this construction is that, the Jacobian of $B_2$
having
trivial 3-torsion,  we reach a curve that was inaccessible before.
We may compose again and again e.g. with parameterizations as in
Section~\ref{sec:61}. It is natural to ask if  we can reach
that way all
genus 2 curves over a large enough finite field or cardinality $q$ when
$q$ is prime to $30$. Answering this question requires to study
some morphisms from a moduli space of covers to
the moduli space of genus 2 curves : proving
in particular that the morphism is surjective
and that the geometric fibers are absolutely irreducible.

\bibliographystyle{plain}

\end{document}